\documentclass{article}
\usepackage{commath}
\usepackage{mathtools}
\usepackage{amsthm}
\usepackage{pgfplots}
\usepackage[round]{natbib}
\usepackage{pgfplotstable}
\usepackage{microtype}
\usepackage[hyphens]{url}
\usepackage{hyperref}

\usepackage[linesnumbered,ruled,vlined]{algorithm2e}
\DontPrintSemicolon

\SetKwComment{Comment}{\color{green!50!black}\# }{}

\newcommand{\assign}{\leftarrow}

\SetKwProg{Function}{def}{}{}

\newcommand{\Reals}{\mathbf{R}}
\newcommand{\Integers}{\mathbf{Z}}
\newcommand{\Naturals}{\mathbf{N}}
\newcommand{\one}{\mathbf{1}}
\newcommand{\Ell}{\mathcal{L}}
\pgfplotsset{compat=1.17}

\mathtoolsset{showonlyrefs}

\newtheorem{theorem}{Theorem}

\theoremstyle{definition}

\newtheorem{problem}{Problem}

\newtheorem{corollary}[theorem]{Corollary}

\theoremstyle{remark}
\newtheorem*{remark}{Remark}

\theoremstyle{remark}
\newtheorem{example}{Example}

\author{Riccardo Sven Risuleo \\ \texttt{riccardosven@gmail.com}}
\title{On the optimal consensus of crab submarines in one dimension}

\begin{document}
\maketitle

\begin{abstract}
    We consider the problem of computing the optimal meeting point of a set of $N$ crab submarines. First, we
    analyze the case where the submarines are allowed any position on the real line: we provide a constructive
    proof of optimality and we use it to provide a linear-time algorithm to find the optimal meeting point for the case
    of sorted starting points.  Second, we use the results for the continuous case to solve the case where the crab
    submarines are restricted to integer locations: we show that, given the solution of the corresponding continuous
    problem, we can find the optimal integer solution in linear time.
\end{abstract}

\section{Introduction}
The crab submarine consensus problem was introduced in the Day 7 of the 2021 \emph{Advent of code} coding
competition~\citep{adventofcode2021day7}. The challenge consists in finding the integer position
where a set of crab submarines can convene starting from different initial positions. The requirement is that the meeting point
should be chosen such that the overall fuel consumption of all crab submarines is minimized.

In the first challenge, the fuel consumption of the crab submarines is linear in the distance traveled, and the solution
(both in the discrete and in the continuous case) is any median of the starting positions.

In the second challenge, the fuel consumption per unit step grows linearly in the distance traveled, so each step costs
1 more fuel than the previous step: the first step costs 1, the second step costs 2, the third step costs 3, and so on.
In the integer case, the total fuel consumption for one crab submarine traveling a distance of $d$ steps is
\begin{equation}
    F = \sum_{s = 0}^d s = \frac{d (d+1)}{2};
\end{equation}
by convexity of the total fuel consumption formula, the integer solution can be found efficiently using, for instance,
bisection search.

As the problem is very new, not many attempts at analytical solutions have been made: as originally formulated, the challenge can
be solved using brute force algorithms (more refined variants based on bisection search are also possible). Notable
exception is the attempt by~\citet{crashandsideburns2021}, who proposed a bound on the distance between the optimal
solution and the algebraic mean of the starting locations. Using a different approach, 
Daniel Szoboszlay\footnote{personal communication.} proposed a relaxation of the problem to the real line with a differentiable cost function.

In this note, we consider the continuous case, where we relax the condition of integer locations and consider crab
submarines that are free to move to any location on the real line. Then, we use the solution to the relaxed problem to
find the optimal meeting point in the integer case.

\section{Problem statement}
A \emph{crab submarine} is a device capable of moving along the real line at arbitrary speed. The crab submarine is
equipped with an engine whose fuel consumption $f$ depends only on the distance traveled $d$ according to
\begin{equation}
    f(d) = \frac{d(d+1)}{2}.
\end{equation}

Consider now a group of $N$ crab submarines, initially at locations $x_i\in \Reals$ for $i=1,...,N$; we want to
determine the meeting location $x^\ast\in \Reals$ such that the overall fuel consumption of all the crab submarines is
minimized. In other words, we want to solve the optimization problem
\begin{equation}
    x^\ast = \arg \min_{x\in\Reals} \Ell(x)
\end{equation}
where
\begin{equation}\label{eq:ell}
    \Ell(x) = \sum_{i=1}^N f\del{\;\envert{x - x_i}\;} = \frac{1}{2} \sum_{i=1}^N (x - x_i)^2 +
    \envert{x-x_i}.
\end{equation}

\begin{remark}
    The function $\Ell(x)$ is continuous, nondifferentiable, and convex.
\end{remark}

\begin{problem}\label{prob:continuous}
    Let $x_1,\ldots, x_N$, $x_i \in \Reals$, be the starting locations of $N$ crab submarines; find $x^\ast \in
    \Reals$ such that $x^\ast = \arg \min_{x\in\Reals} \Ell(x)$.
\end{problem}

Restricting the function $\Ell(x)$ to the integers, we obtain the following discrete crab submarine problem:

\begin{problem}\label{prob:discrete}
    Let $k_1,\ldots, k_N$, $k_i \in \Integers$, be the starting locations of $N$ crab submarines; find $k^\ast \in
    \Integers$ such that $k^\ast = \arg \min_{k \in \Integers} \Ell(k)$.
\end{problem}

\section{Prerequisites}
To solve the problems, we need some definitions. First, we define the \emph{mean} of the crab submarine starting
positions as the point $\tilde x\in \Reals$ such that
\begin{equation}
    \tilde x = \frac{1}{N}\sum_{i=1}^N x_i.
\end{equation}
Furthermore a \emph{median} of the crab positions as any point $\hat x \in \Reals$ such that
\begin{equation}\label{eq:median}
    \sum_{i=1}^N \one_{\sbr{x_i \leq \hat x}} = \sum_{i=1}^N \one_{\sbr{x_i \geq \hat x}}
\end{equation}
where $\one_A$ is the indicator variable of the proposition $A$:
\begin{equation}
    \one_A = \begin{cases}
        1 & \text{if } A, \\
        0 & \text{otherwise}.
    \end{cases}
\end{equation}

\begin{remark}
    If $N$ is odd, there exists an index $m$ such that $\hat x = x_m$; this coincides with the standard definition of
    the median of a set of real numbers. If $N$ is even, there exists two indices $\underline{m}$ and $\overline{m}$
    such that any $x \in \sbr{x_{\underline m},\, x_{\overline m}}$ is a median.
\end{remark}

In addition, we define the count imbalance $\Delta(x)$ as the difference between the number of crab starting locations
strictly to the left and strictly to the right of point $x$:
\begin{equation}
    \Delta(x) = \sum_{i=1}^N \one_{\sbr{x_i < x}} - \sum_{i=1}^N \one_{\sbr{x_i > x}}.
\end{equation}
with this definition, a median $\hat x$ is any point such that $\Delta(\hat x) = 0$. We also define the count
multiplicity as the number of crab starting locations at a point $x$:
\begin{equation}
    \rho(x) = \sum_{i=1}^N \one_{\sbr{x_i = x}}.
\end{equation}

Finally, a \emph{subgradient} is a generalization of the gradient for non differentiable
functions~\citep{shor2012minimization}. The subgradient
$\partial g(x)$ of a function $g(x)$ is the set of all the slopes of the tangents to the function at $x$; in other
words, it is the set of real numbers $c$ such that
\begin{equation}
    g(y) - g(x) \geq c (y - x)
\end{equation}
for all $y$. If $g$ is \emph{convex}, the subgradient is a nonempty closed interval $\sbr{a,\, b}$, where $a$ and $b$
are given by
\begin{equation}
    a = \lim_{y \to x^-} \frac{f(y) - f(x)}{y - x},\qquad b = \lim_{y \to x^+} \frac{f(y) - f(x)}{y - x}.
\end{equation}
The subgradient over a domain is a singleton set if and only if the function is differentiable over that domain.

For example, the function $\envert{x}$ is not differentiable over $\Reals$. However, it has a subgradient given by
\begin{equation}
    \partial\envert{x} = \begin{cases}
        \{-1\} & \text{if } x < 0, \\
        \{1\} & \text{if } x > 0, \\
        \sbr{-1, \, 1} & \text{if } x = 0.
    \end{cases}
\end{equation}

Using subgradients, we can extend the necessary first-order conditions for optimality to state that a point $x^\star$ is
the minimum of a convex function $g$ if 0 belongs to the subgradient of $g$ at
$x^\star$~\citep[see][]{rockafellar2015convex}:
\begin{equation}
    0 \in \partial g(x^\star).
\end{equation}

To proceed with the solution, we first compute the subgradient of~\eqref{eq:ell}
\begin{align}
    \partial \Ell(x) &= \frac{1}{2}\sum_{i=1}^N 2(x - x_i) + \partial\envert{x - x_i} \\
                     &= N x - N \tilde x + \frac{1}{2}\sum_{i=1}^N\partial\envert{x - x_i} \\
    &= \begin{cases}
        N (x - \tilde x) + \frac{\Delta(x)}{2} + \sbr[1]{-\frac{\rho(x)}{2},\,\frac{\rho(x)}{2}}& \text{if } x\in\cbr{x_1,x_2,\ldots,
        x_N}, \\
            \cbr{N(x - \tilde x) + \frac{\Delta(x)}{2}}  & \text{otherwise}.
    \end{cases}
\end{align}
Note that, in these expressions, ``$+$'' denotes Minkowski sum where appropriate.

\section{Results}
We now present the main results. We first focus on the continuous case and we present a constructive theorem for the
optimal solution to Problem~\ref{prob:continuous}; then, we provide a constructive solution to
Problem~\ref{prob:discrete}.

\begin{theorem}\label{thm:solution}
    If there exists an index $t$ such that
    \begin{equation}\label{eq:condition-t}
        \envert{N(x_t - \tilde x) + \frac{\Delta(x_t)}{2}} \leq \frac{\rho(x_t)}{2},
    \end{equation}
    then $x^\ast = x_t$ is the solution to Problem~\ref{prob:continuous}. If no such index exists, then $x^\ast$ is the solution to
    the following equation:
    \begin{equation}\label{eq:condition-ast}
        x^\ast + \frac{\Delta(x^\ast)}{2N} = \tilde x .
    \end{equation}
\end{theorem}
\begin{proof}
    Consider the subgradient $\partial \Ell(x)$ and consider $x \in \cbr{x_1,\ldots,x_N}$; for all these points, we have that
    \begin{equation}
        \partial \Ell(x) = \sbr{N(x - \tilde x) + \frac{\Delta(x) - \rho(x)}{2},\; N(x - \tilde x)
        +\frac{\Delta(x) + \rho(x)}{2}}.
    \end{equation}
    Suppose now that there exists an index $t$ such that~\eqref{eq:condition-t} holds, then we have that
    \begin{equation}
        N(x_t - \tilde x) + \frac{\Delta(x_t)-\rho(x_t)}{2} \leq 0, \qquad 
        N(x_t - \tilde x) + \frac{\Delta(x_t)+\rho(x_t)}{2} \geq 0,
    \end{equation}
    hence, $0 \in \partial\Ell(x_t)$ and we have the proof. For any $x \notin \cbr{x_1, \ldots, x_N}$, the subgradient
    contains only one point, and $0 \in \partial \Ell(x^\ast)$ if and only if there exists a solution $x^\ast$
    to~\eqref{eq:condition-ast}.
\end{proof}

\begin{corollary}
    Let $x^\ast$ be a solution to Problem~\ref{prob:continuous}, then
    \begin{equation}
        \envert{N(x^\ast - \tilde x) + \frac{\Delta(x^\ast)}{2}} \leq \frac{\rho(x^\ast)}{2},
    \end{equation}
\end{corollary}
\begin{proof}
    Follows from the fact that $\rho(x) = 0$ for $x\notin\cbr{x_1, \ldots, x_N}$.
\end{proof}

\begin{corollary}\label{cor:mean}
    If the mean starting position $\tilde x$ is a median, then $x^\ast = \tilde x$ is a solution to
    Problem~\ref{prob:continuous}.
\end{corollary}
\begin{proof}
    If $\tilde x$ is a median, then $\Delta(\tilde x)= 0$. Suppose that $x_t = \tilde x$ for some $t$; then $\rho(\tilde x)> 0$
    and~\eqref{eq:condition-t} is verified; otherwise,~\eqref{eq:condition-ast} reduces to $x^\ast = \tilde x$.
\end{proof}

\begin{corollary}\label{cor:median}
    If there exists an index $m$ such that $x_m$ is a median and
    \begin{equation}\label{eq:median-thm}
        \envert{x_m - \tilde x} \leq \frac{\rho(x_m)}{2N},
    \end{equation}
    then $x^\ast = x_m$ is a solution to Problem~\eqref{prob:continuous}.
\end{corollary}
\begin{proof}
    Follows from the same argument as the proof of Corollary~\ref{cor:mean} applied to~\eqref{eq:condition-t}.
\end{proof}

Note that, using the convexity of $\Ell(x)$ and Theorem~\ref{thm:solution}, it is possible to solve the discrete crab-submarine case:
\begin{theorem}\label{thm:discrete}
    Let $x^\ast \in \Reals$ be a solution to Problem~\ref{prob:continuous} and let $k^-  = \lfloor x^\ast \rfloor \in
    \Integers$ be the largest integer smaller than or equal to
    $x^\ast$; similarly, let $k^+ = \lceil x^\ast\rceil \in \Integers$ be the smallest integer larger than or equal to
    $x^\ast$. Then,
    \begin{equation}
        k^\ast = \arg \min_{k \in \cbr{k^-,\, k^+}}\Ell(k),
    \end{equation}
    is the solution to Problem~\ref{prob:discrete}.
\end{theorem}
\begin{proof}
    We consider the case that $x^\ast \notin \cbr{\lfloor x^\ast \rfloor,\, \lceil x^\ast\rceil }$: the proof is
    trivial otherwise.

    Let $n\in \Naturals$, and let $k^+ = \lceil x^\ast \rceil$; then,we have that $x^\ast < k^+ < k^+ + n$. By the
    definition of convexity, for all $t\in \sbr{0,\,1}$,
    \begin{align}
        \Ell\del{ t x^\ast + (1 - t)(k^+ + n)} &\leq t \Ell(x^\ast) + (1-t)\Ell(k^+ + n) \\ 
                                               & \leq \Ell(k^+ + n)
    \end{align}
    where we have used the fact that $\Ell(x^\ast) \leq \Ell(k^+ + n)$. Then, let
    \begin{equation}
        t = \frac{n}{k^+ + n - x^\ast},
    \end{equation}
    we have that $0 < t < 1$; substituting into the previous inequality, we have
    \begin{equation}
        \Ell(k^+) \leq \Ell(k^+ +n),
    \end{equation}
    for all $n > 0$; which proves that $k^+$ is the minimum of $\Ell(k)$ for all $k\in \Integers,\, k > x^\ast$. A
    similar argument applied to $k^- - n < k^- < x^\ast$, shows that $k^-$ is the minimum of $\Ell(k)$ for all $k
    \in \Integers,\, k < x^\ast$. Hence, the minimum of $\Ell(k)$ must be in $\cbr{k^-, \,k^+}$ and we have the
    proof.
\end{proof}

Theorem~\ref{thm:discrete} shows that we can find the solution to Problem~\ref{prob:discrete} in linear time, once we
have the solution to Problem~\ref{thm:solution}, by checking the cost function $\Ell(x)$ in the closest integers on
each side of $x^\ast$. As evaluating $\Ell(x)$ is $O(N)$, we can solve Problem~\ref{prob:discrete} in linear time if
we have $x^\ast$; in Section~\ref{sec:linear}, we present a linear-time algorithm to solve
Problem~\ref{prob:continuous} for the case of sorted starting locations.

% However, the convexity of the cost function allows us to give a stronger result:
% \begin{theorem}\label{thm:discrete}
%     Let $x^\ast \in \Reals$ be a solution to Problem~\ref{prob:continuous}. If $x^\ast \in \Integers$, then $k^\ast =
%     x^\ast$ is a solution to Problem~\ref{prob:discrete}; if not, the solution is
%     \begin{equation}
%         k^\ast = \begin{cases}
%             k^+ & \text{if } x^\ast - \tilde x -\frac{\Delta(x^\ast)}{2N} \geq 0, \\ 
%             k^- & \text{otherwise}.
%         \end{cases}
%     \end{equation}
% \end{theorem}
% \begin{proof}
%     If $x^\ast\in\Integers$, then  $k^+ = k^- = x^\ast$ and we have the proof from Lemma~\ref{lem:discrete}.
% \end{proof}
% 
% 
\section{Examples}
The first example shows a case where the mean is also a median and is the optimal solution.
\begin{example}\label{ex:one}
    Suppose that two crab submarines start at locations $x_1 = 0$ and $x_2 = 1$, then $\tilde x = 0.5$ is
    a median and is also the optimal meeting point according to Corollary~\ref{cor:mean}.
\end{example}

The second example shows a case where a starting location is a median, and is a solution.
\begin{example}\label{ex:two}
    Suppose that three crab submarines start at locations $x_1 = 0$, $x_2 = 1/3$, and $x_3 = 1/2$, then $\tilde x =
    5/18$, and $\Delta(x_2) =  0$: in this case $x_2$ is a median, furthermore
    \begin{equation}
        \envert{x_2 - \tilde x} = \envert{\frac{1}{3} - \frac{5}{18}} = \frac{1}{18} < \frac{1}{2N},
    \end{equation}
    so $x_2$ is an optimal meeting point according to Corollary~\ref{cor:median}.
\end{example}

The third example shows a case where none of the starting locations is a solution.
\begin{example}\label{ex:three}
    Suppose that three crab submarines start at location $x_1 = 0$, $x_2 = 1$, and $x_3 = 1$, then $\tilde x = 2/3$,
    then $\Delta(x_1) = -2$, $\Delta(x_2) = \Delta(x_3) = 1$; we have
    \begin{equation}
        \envert{x_1 - \tilde x + \frac{\Delta(x_1)}{2N}} = \envert{0 - \frac{2}{3} - \frac{2}{6}} = 1 >
        \frac{1}{6},
    \end{equation}
    so $x_1$ is not a solution according to~\ref{thm:solution}; however:
    \begin{equation}
        \envert{x_2 - \tilde x + \frac{\Delta(x_2)}{2N}} = 
        \envert{x_3 - \tilde x + \frac{\Delta(x_3)}{2N}} = \envert{1 - \frac{2}{3} + \frac{1}{6}} = \frac{1}{2} >
        \frac{1}{2N},
    \end{equation}
    so neither $x_2 = x_3$ are solutions. For this problem, we have that
    \begin{equation}
        \Delta(x) = \begin{cases}
            -3 & x < 0 \\
            -2 & x = 0 \\
            -1 & 0 < x < 1 \\
            +3 & x > 1
        \end{cases}
    \end{equation}
    Therefore, it can be seen that $x^\ast = 5/6$ solves~\eqref{eq:condition-ast} and is the solution to
    Problem~\eqref{prob:continuous}.
\end{example}

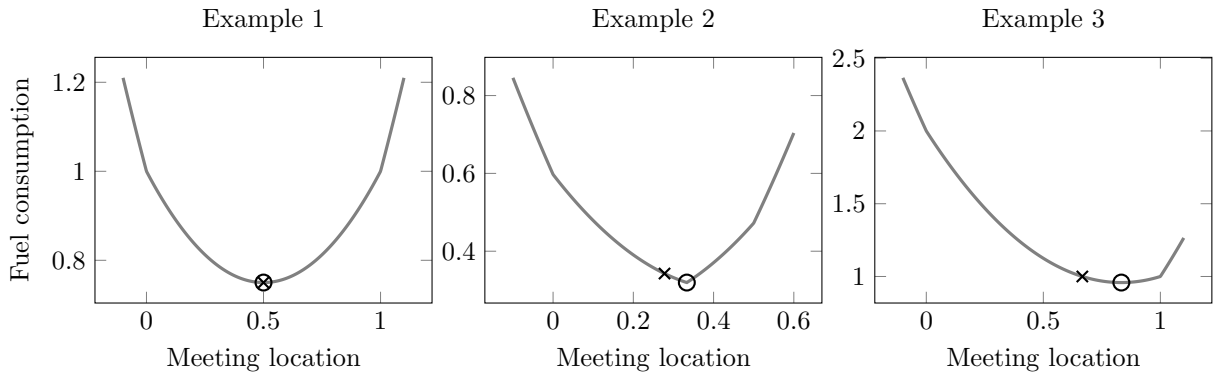
\begin{figure}[tb]
    \centerline{
    \begin{tikzpicture}
        \pgfplotstableread[col sep=comma]{simulations/example1.csv}{\datatableA};
        \pgfplotstableread[col sep=comma]{simulations/example2.csv}{\datatableB};
        \pgfplotstableread[col sep=comma]{simulations/example3.csv}{\datatableC};
        \begin{axis}
            [ylabel=Fuel consumption,
                xlabel=Meeting location,
                width=0.5\textwidth,
                height=0.4\textwidth,
                name=A,
                ylabel near ticks,
                title=Example 1,
            ]
            \addplot[solid, gray, very thick] table[x=x, y=y] {\datatableA};
            \addplot[only marks, thick, mark=o, mark size=3pt] coordinates { (0.5, 0.75)};
            \addplot[only marks, thick, mark=x, mark size=3pt] coordinates { (0.5, 0.75)};
        \end{axis}
        \begin{axis}
            [
                xlabel=Meeting location,
                width=0.5\textwidth,
                height=0.4\textwidth,
                name=B,
                anchor=north west,
                at=(A.north east),
                xshift = 2em,
                title=Example 2,
            ]
            \addplot[solid, gray, very thick] table[x=x, y=y] {\datatableB};
            \addplot[only marks, thick, mark=o, mark size=3pt] coordinates { (0.3333333333333333, 0.3194444444444444)};
            \addplot[only marks, thick, mark=x, mark size=3pt] coordinates { (0.27777777777777773, 0.34259259259259256)};
        \end{axis}
        \begin{axis}
            [
                xlabel=Meeting location,
                width=0.5\textwidth,
                height=0.4\textwidth,
                anchor = north west,
                at=(B.north east),
                xshift = 2em,
                title=Example 3,
            ]
            \addplot[solid, gray, very thick] table[x=x, y=y] {\datatableC};
            \addplot[only marks,thick, mark=o, mark size=3pt] coordinates { (0.8333333333333334, 0.9583333333333334)};
            \addplot[only marks,thick, mark=x, mark size=3pt] coordinates { ( 0.6666666666666666,1.0000000000000002)};
        \end{axis}
    \end{tikzpicture}
}
    \caption{Simulation results for the examples: the plots show fuel consumption as a function of the meeting point.
    The cross shows the mean starting location $\tilde x$; the circle shows the optimal point.
$x^\ast$.}\label{fig:simulation}
\end{figure}

Simulation results showing the function $\Ell(x)$, together with the mean starting locations and the optimal meeting
points are presented in Figure~\ref{fig:simulation}.

\section{A linear-time algorithm}\label{sec:linear}
Using Theorem~\ref{thm:solution}, and the assumption that the starting locations are indexed such that $x_i \leq
x_{i+1}$ for all indices $i$, we can define a linear-time algorithm to find the optimal meeting point for $N$ crab
submarines on the real line. A pseudocode implementation is presented in Algorithm~\ref{alg:pseudocode}. Further
reduction of the computational complexity (at the expense of a linear space complexity) is possible, by pre-computing
$\Delta(x)$ and $\rho(x)$ for all starting locations and replacing the linear search with a bisection search.

\begin{algorithm}
    \caption{Compute the optimal meeting location of $N$ crab submarines; convexity of the function $\Ell(x)$ guarantees
    that there is a solution.}\label{alg:pseudocode}
    \Function{CrabSubmarines($x_1,\ldots, x_N$)}{
        $\tilde x \assign \frac{1}{N}\sum_{i=1}^N x_i$ \Comment{Mean starting location}
        $\rho \assign 1$ \Comment{Counts for repeated starting locations}
        %$x_{N+1} \assign \infty$\;
        $d \assign -1/2$ \Comment{Keeps track of $\frac{\Delta(x)}{2N}$}
        \For{$i = 1, \ldots, N$}{
            \If(\Comment*[h]{$x_{N+1}:=\infty$}){$x_i = x_{i+1}$}{
            $\rho \assign \rho + 1$
        }\Else{
               $d \gets d + \frac{\rho}{2N}$ \;
               \If(\Comment*[h]{Check condition~\eqref{eq:condition-t}}){$\envert{x_i - \tilde x + d}
                   \leq \frac{\rho}{2N}$}{
                   \Return{$x_i$}
               }
               $d \gets d + \frac{\rho}{2N}$ \;
               $x \gets \tilde x - d$ \Comment{Candidate solution to~\eqref{eq:condition-ast}}
               \If(\Comment*[h]{$x_{N+1}:=\infty$}){$x_i \leq x \leq x_{i+1}$}{
                   \Return{$x$}
               }
               $\rho \gets 1$
           }
       }
    }
\end{algorithm}

\begin{remark}
    The running cost of the proposed algorithm is $O(N)$; however, this cost is conditioned on the fact that the starting
    locations $x_1,\ldots,x_N$ are ordered. If that is not the case, the input can be sorted before running
    Algorithm~\ref{alg:pseudocode}, for an asymptotic complexity of $O(N\log N)$.
\end{remark}

\section{Conclusions}
We have established necessary and sufficient conditions for optimality for the meeting point of an arbitrary finite
number of crab submarines and, in the case where the starting locations are ordered, we have proposed a linear-time
algorithm for computing the optimal meeting point on the real line. Also, we have shown how the solution to the
relaxation of the discrete problem to the real line can be used to find the solution to the discrete problem in linear
time.

We have various interesting points of future research:
\begin{itemize}
    \item Can we extend the results for an infinite number of crab submarines? Suppose that we have a continuous
        density $\mu(x)$ of crab submarines, how is the optimal meeting point related to the mean, the mode, and the
        median of the density?
    \item Can we extend the results for crab submarines in multiple dimensions? Do linear-time algorithms still exist?
\end{itemize}

\section*{Acknowledgements}
We thank an anonymous reviewer for their vital input on multiple parts of the first draft.

\bibliographystyle{plainnat}
\bibliography{main}

\end{document}